\input amstex
\input amsppt.sty
\magnification\magstep1

\def\ni\noindent
\def\sbs{\subset}

\def\dim{\operatorname{dim}}

\def\R{\text{\bf R}}

\def\Z{\text{\bf Z}}

\def\sU{\Cal U}

\def\sZ{\Cal Z}

\hoffset= 0.0in
\voffset= 0.0in
\hsize=32pc
\vsize=38pc
\baselineskip=24pt
\NoBlackBoxes
\topmatter
\author
A.N. Dranishnikov
\endauthor

\title
On Bestvina-Mess formula 
\endtitle
\abstract

Bestvina and Mess [BM] proved a remarkable formula for  torsion free
hyperbolic groups
$$
\dim_L\partial\Gamma=cd_L\Gamma-1
$$
connecting the cohomological dimension of a group $\Gamma$ with
the cohomological dimension of its boundary $\partial\Gamma$.
In [Be] Bestvina introduced a notion of
$\sZ$-structure on a discrete group and noticed that his formula holds true for all
torsion free groups with $\sZ$-structure.
Bestvina's notion of $\sZ$-structure can be extended to groups containing
torsion by replacing the covering space action in the
definition by the geometric action. Though the Bestvina-Mess
formula trivially is not valid for
groups with torsion,
we show that it still holds in the following modified form:
{\it The cohomological dimension of a $\sZ$-boundary of a group $\Gamma$
equals
its global cohomological dimension for every PID $L$ as the coefficient group}
$$
\dim_L\partial\Gamma=gcd_L(\partial\Gamma).
$$
Using this formula we show that the cohomological dimension of
the boundary $\dim_{L}\partial\Gamma$ is a quasi-isometry invariant
of a group.

For CAT(0) groups and $L=\Z$ these results were obtained in [GO].
\endabstract

\thanks The author was partially supported by NSF grants DMS-0305152
\endthanks

\address University of Florida, Department of Mathematics, P.O.~Box~118105,
358 Little Hall, Gainesville, FL 32611-8105, USA
\endaddress


\email  dranish\@math.ufl.edu
\endemail

\endtopmatter

\document
\head \S1 Introduction \endhead
A closed subset  $Z$ of an $AR$-space $\bar X$ is called a {\it Z-set}
if there is a homotopy $H:\bar X\times[0,1]\to\bar X$ such that
$H_0$ is the identity and $H_t(\bar X)\subset\bar X\setminus Z$ for all
$t>0$.
An equivalent statement is that
$Z$ is a Z-set if for every open set $U\subset\bar X$ the inclusion
$U\setminus Z\subset U$ is a homotopy equivalence.

The following generalizes Bestvina's definition from \cite{Be}.

DEFINITION. A discrete group $\Gamma$ has $\sZ$-structure if there is a pair
$(\bar X,Z)$ of compact spaces satisfying the axioms:
\roster
\item{} $\bar X$ is AR;
\item{} $Z$ is a Z-set in $\bar X$;
\item{} $X=\bar X\setminus Z$ is a metric space that
admits a proper discontinuous action of $\Gamma$  by isometries with
a compact quotient;
\item{} The compactification $\bar X$ of $X$ is Higson dominated.
\endroster

Originally Bestvina required in (3) a covering space action which
automatically makes $\Gamma$ torsion free.
Note that for torsion free groups this definition of $\sZ$-structure
coincides with the original one.

The {\it Higson compactification} $\bar X_H$ [Roe] of a metric
space $X$ is defined by means of the ring $C_H(X)$ of bounded
continuous functions $f:X\to\R$ with variation tending to zero at
infinity. Namely, given $R<\infty$ and $\epsilon>0$, there is a
compact set $K\subset X$ such that $diam(f(B_R(x))<\epsilon$ for
$x\in X\setminus K$ where $B_R(x)$ denotes the $R$-ball centered
at $x$. A compactification $\bar X$ of $X$ is called {\it Higson
dominated} if the identity map $id_X$ extends to a continuous map
$\alpha:\bar X_H\to\bar X$.

For a group $\Gamma$ with Z-structure $(\bar X,Z)$ we will denote the set
$Z$ as $\partial\Gamma$ and call it a {\it boundary of the group} $\Gamma$.

For a set $A\subset X$ the trace $A'$ of $A$ on the boundary $Z$
is the intersection with $Z$ of the closure of $A$ in $\bar X$,
 $$A'=Cl_{\bar X}A\cap Z.$$
The Axiom 4 can be restated as follows:

(4) {\it For every set $A\subset X$ and every $r>0$ the trace of $A$ in $Z$
coincides with the trace of the $r$-neighborhood $N_r(A)$}:
$$
A'=(N_r(A))'.
$$
We note that the condition (4) is equivalent to the condition that
the action of $\Gamma$ is
{\it small at infinity} [DF] i.e.,
for every $x\in Z$ and a neighborhood $U$ of $x$ in $\bar X$,
for every compact set $C\subset X$ there is a smaller neighborhood $V$,
$x\in V\subset U$, such that
$g(C)\cap V\ne\emptyset$ implies $g(C)\subset U$ for all $g\in\Gamma$.
So we can reformulate Axiom 4 once more:

(4) {\it For every $x\in X$ and every $r>0$ the sequence $\{B_r(\gamma(x))\mid
\gamma\in\Gamma\}$ is a 0-sequence in $\bar X$.}

The Bestvina-Mess formula proven first for torsion free hyperbolic groups
[BM] as it noted in [Be] is valid for all torsion free groups with
$\sZ$-structure.
It states
$$
\dim_LZ=cd_L\Gamma-1
$$
where $L$ ia a PID, $dim_LZ$ is the cohomological dimension of a space $Z$ with
coefficients in $L$, and $cd_L\Gamma$ is the cohomological dimension of a group
$\Gamma$ with coefficient in $L$. For a torsion free group $\Gamma$ with a
$\sZ$-structure $(\bar X,Z)$ we have [Be]
$$
cd_L\Gamma=\max\{n\mid H^n_c(X;L)\ne 0\}.
$$

We recall that a {\it global cohomological dimension} of a space
$X$ with coefficients in an abelian group $G$ is the following
number
$$
gcd_G(X)=\max\{n\mid H^n(X;G)\ne 0\}.
$$
Since $\bar X$ is contractible, the exact sequence of the pair $(\bar X,Z)$
implies that
$$
cd_L\Gamma-1=gcd_L(Z).
$$
Thus, the Bestvina-Mess formula is equivalent to the equality
$$
dim_LZ=gcd_L(Z).
$$
The purpose of this paper is to prove this equality for general groups $\Gamma$
with $\sZ$-structure. We do it by an adjustment of Bestvina-Mess' argument.
In the case when $\Gamma$ is a CAT(0) group and $L=\Z$ this formula was
established by Geoghegan and Ontaneda [GO] by a different method.

\head  \S2 Global and local cohomological dimension of a $\sZ$-boundary\endhead

The main result is the following 
\proclaim{Theorem 1} Let $(\bar
X,Z)$ be a $\sZ$-structure on a group $\Gamma$. Then for every
principle ideal domain $R$ the cohomological dimension of $Z$
coincides with its global cohomological dimension:
$$
dim_RZ=gcd_R(Z).
$$
\endproclaim

We recall that the {\it cohomological dimension} with respect to the coefficient group
$G$ of a locally compact metric space $Z$ is the following number
$$
dim_GZ=\max\{n\mid H^n_c(U;G)\ne 0, U\subset_{Open} Z\}
$$
where $H^*_c$ denotes the cohomology with compact supports.

Let $\dim_GX=r$. We say that a point $x\in X$ is {\it
dimensionally essential} (with respect to $G$) if there is a
neighborhood $U$ of $x$ such that for every smaller neighborhood
$V$ of $x$ the image of the inclusion homomorphism
$i_{V,U}:H^r_c(V;G)\to H^r_c(U;G)$ is nonzero. Alexandroff called
such points as containing an $r$-dimensional obstruction. We refer
to [Dr] for more details on cohomological dimension as well for
the following Lemma. 
\proclaim{Lemma 1}(P.S. Alexandroff) Let $X$
be a compact metric space with $\dim_GX=r$. Then the set of
dimensionally essential points in $X$ contains a locally compact
subset $Y$ of $\dim_GY=r$.
\endproclaim
\demo{Proof}
Let $W$ be an open subset of $X$ with $H_c^r(W;G)\neq 0$. Because
of the continuity of cohomology there is a closed in $W$
set $Y$ minimal with respect the property: The inclusion
homomorphism $H^r_c(W;G)\to H^r_c(Y;G)$ is nonzero. Then, clearly, 
$\dim_GY=r$.  For every $x\in Y$ we
take $U=W$. Let $V\subset U$ be a neighborhood of $x$.
Consider the diagram generated by
exact sequences of pairs $(U,U\setminus V)$ and $(Y,Y\setminus V)$.
$$
\CD
 H^r_c(V;G) @>{i_{V,U}}>> H_c^r(U;G) @.\\
@Vj_{V,V\cap Y}VV   @Vj_{U,Y}VV\ @.\\
H^r_c(V\cap Y;G)  @>{i_{V\cap Y,Y}}>> H^r_c(Y;G)
@>{j_{Y,Y\setminus V}}>> H^r_c(Y\setminus V;G)\\
\endCD
$$

Let $\alpha\in H^r_c(U;G)$ such that $j_{U,Y}(\alpha)\neq 0$.
Since $Y$ is minimal, $j_{Y,Y\setminus V}(j_{U,Y}(\alpha))=0$.
The exactness of the bottom row implies that there is
$\beta\in H^r_c(Y\cap V;G)$ such that $i_{Y\cap V,Y}(\beta)=
j_{U,Y}(\alpha)$. Since $\dim_GV\leq r$, the homomorphism
$j_{V,Y\cap V}$ is an epimorphism and hence there is
$\gamma\in H^r_c(V;G)$ with $j_{V,Y\cap V}(\gamma)=\beta$.
Therefore $j_{U,Y}i_{V,U}(\gamma)\neq 0$ and hence
$i_{V,U}(\gamma)\neq 0$. \qed
\enddemo
The following lemma is contained essentially in [BM] (Proposition 2.6).
\proclaim{Lemma 2}
Suppose that $Z$ is a Z-set in an AR space $\bar X$ and $\dim_RZ=n$
where $R$ is an abelian group.
Then for every dimensionally essential point $z\in Z$
and for every neighborhood $W\subset \bar X$
there is an open neighborhood $\bar U\subset\bar X$, $z\in \bar U\subset W$,
such that the coboundary homomorphism for the pair $(\bar U,\bar U\cap Z)$
$$
H^n_c(V;R)\to H_c^{n+1}(U;R)
$$
is nonzero monomorphism where $V=\bar U\cap Z$ and $U=\bar U\setminus V$.
\endproclaim
\demo{Proof}
Given $W$ we take open in $Z$ set $V$, $Cl(V)\subset W$,
with $H^n_c(V;R)\ne 0$. Let $A=Z\setminus V$.
In view of the $Z$-set condition there is a homotopy
$F:\bar X\times [0,1]\to \bar X$ such that
\roster
\item{} $F_0$ is the identity;
\item{} $F_t|_A$ is the inclusion $A\subset Z$;
\item{} $F_t(\bar X\setminus A)\subset \bar X\setminus Z$ for all $t>0$;
\item{} $\bar X\setminus W\subset F_1(\bar X)$.
\endroster
Let $Y_1=F_1(\bar X)$ and $Y_2=F(Y_1\times[0,1])$. Then
$Y_1\subset Y_2$, $Y_1\cap Z= Y_2\cap Z=A$. Since the inclusion
$(\bar X,Y_1)\to (\bar X,Y_2)$ can be homotoped to the map to
$(Y_2,Y_2)$, the homomorphism $H^n(\bar X,Y_2;R)\to H^n(\bar
X,Y_1;R)$ is zero. Consider the diagram
$$
\CD
H^n(\bar X,Y_2;R) @>>> H^n(Z\cup Y_2,Y_2;R) @>j>> H^{n+1}(\bar X,Z\cup Y_2;R)\\
@V{0}VV @V{\cong}VV @VVV\\
H^n(\bar X,Y_1;R) @>>> H^n(Z\cup Y_1,Y_1;R) @>>> H^{n+1}(\bar X,Z\cup Y_1;R)\\
\endCD
$$
The diagram chasing shows that $j$ is a monomorphism.
We define $\bar U=\bar X\setminus Y_2$. Then $U=\bar U\setminus Z$, $V=\bar U\cap Z$, $H^n(Z\cup Y_2,Y_2;R)=H^n_c(V;R)$,
$H^{n+1}(\bar X, Z\cup Y_2)=H_c^{n+1}(U;R)$, and $j$ is
the coboundary homomorphism from the exact sequence of the pair $(\bar U,V)$.
Since $j$ is a monomorphism and  $H^{n}_c(V;R)\ne 0$,
it is nonzero.
\qed
\enddemo
We recall that a metric space $X$ is {\it uniformly contractible}
if there is a function $\rho:\R_+\to\R_+$ with the property that
every ball $B_r(x)$ of radius $r$ centered at $x$ is contractible
inside the ball $B_{\rho(r)}(x)$. A metric space $X$ has {\it
bounded geometry} if there is $\epsilon>0$ such that for every
$R<\infty$ the $\epsilon$-capacity of $R$-balls $B_R(x)$ is
uniformly bounded from above. The following is obvious.
\proclaim{Proposition 1} Let $\Gamma$ act properly discontinuously
and cocompactly by isometries on a proper metric AR-space $X$ then
$X$ is uniformly contractible with bounded geometry.
\endproclaim
\proclaim{Corollary 1}
Let $X$ as in Proposition 1, then $X$ admits a uniformly bounded open cover
$\sU$ with finite dimensional nerve $K=N(\sU)$ and a left proper homotopy
inverse $s:K\to H$ to the projection to the nerve $\psi:X\to K$ with  bounded
homotopy $H:X\times I\to X$ joining $s\circ\psi$ with the identity $id_X$.
\endproclaim
\demo{Proof}
This is a standard fact of a coarse geometry and it can be found in [HR].
\enddemo

We denote by $\dot B_r(x)$ an open ball around $x$ of radius $r$.

DEFINITION.[GO] Let $X$ be a locally compact metric space.
The group $H^i_c(X;G)$ is {\it uniformly trivial} if there is a
function $\nu:\R_+\to\R_+$ such that the inclusion homomorphism
$$
H^i_c(\dot{B}_r(x);G)\to H^i_c(\dot{B}_{\nu(r)}(x);G)
$$
is trivial for all $x\in X$ and $r>0$.

\proclaim{Proposition 2} Let $\Gamma$ act properly discontinuously
and cocompactly by isometries on a proper metric AR-space $X$ and
let $H_c^i(X;R)=0$ for a PID $R$. Then $H_c^i(X;R)=0$ is uniformly
trivial.
\endproclaim
\demo{Proof}
In view of cocompactness of the action, the orbit $\Gamma x_0$ is an
$r_0$-dense in $X$. Hence it suffices to prove that for every $r$ there is
$\nu>0$ such that the homomorphism
$$
H^i_c(\dot{B}_r(x_0);R)\to H^i_c(\dot{B}_{\nu}(x_0);R)
$$
is trivial.
Let $\psi:X\to K$ and $s:K\to X$ be from Corollary 1 and let
$d$ be an upper bound on a homotopy between $s\circ\psi$ and $id_X$.
For a subset $A\subset K$ we denote by $st(A)$ the union of all simplices that
have nonempty intersection with $A$. Let $L=st(\psi(X\setminus B_{r+2d}(x_0)))$.
Then the inclusion
$X\setminus \dot B_{r+2d}(x_0)\to X\setminus \dot B_r(x_0)$ is homotopy factored through
the maps
$$
X\setminus \dot B_{r+2d}(x_0) @>\psi>> L
@>s>>X\setminus \dot B_r(x_0).
$$
Therefore the inclusion homomorphism
$H^i_c(\dot{B}_r(x_0);R)\to H^i_c(\dot{B}_{r+2d}(x_0);R)$
is factored through the group $H^i(K,L;R)$ which is a finitely generated
$R$-module. Let $a_1,\dots a_m\in H^i(K,L;R)$ be a set of generators.
For every $a_j$ there is $\nu_j$ such that the element $\psi^*(a_j)$
vanishes in $H^i_c(B_{\nu_j}(x_0);R)$. We take $\nu=\max\{\nu_j\}$.\qed
\enddemo
The following Lemma is an extension of a claim that appeared in
the proof of Corollary 1.4.(a) [BM]. 
\proclaim{Lemma 3} Let
$\Gamma$ act properly discontinuously and cocompactly by
isometries on a proper metric AR-space $X$ and let $R$ be a PID.
Suppose that $H^i_c(X;R)=0$ for all $i>n$. Then there is a number
$r$ such that for every open set $U\subset X$ the inclusion 
$U\subset N_{r}(U)$ into the
open $r$-neighborhood induces the zero
homomorphism $H_c^i(U;R)\to H_c^i(N_r(U);R)$  for all $i>n$.
\endproclaim
\demo{Proof} We consider a sequence of locally finite uniformly
bounded open coverings of $X$, $\sU_0\succ \sU_1 \succ  \sU_2\succ
\sU_3\ \dots$ where $\sU_{k+1}$ is a refinement of $\sU_k$, with
$\lim_{k\to 0}mesh(\sU_k)=0$, and $\sU_0$ as in Corollary 1. Let
$\psi_j:X\to K_j$ denote a projection to the nerve of $\sU_j$. Let
$s:K_0\to X$ be a proper homotopy lift and let $d>0$ be the size
of a homotopy between $s\circ\psi_0$ and $id_X$.

Let $U\subset X$ be an open set and let $L=st(\psi_0(X\setminus N_{2d}(U)))$.
Then the homotopy commutative diagram
$$
\CD
X\setminus N_{2d}(U) @>>> X\setminus\psi_0^{-1}(L) @>>> X\setminus U\\
@ V{\psi_0}VV @ V{\psi_0}VV @ A{s}AA\\
L @>{=}>> L @>{=}>> L\\
\endCD
$$
defines the commutative diagram for cohomology
$$
\CD
H^i_c(N_{2d}(U);R) @<<< H^i_c(\psi_0^{-1}(K_0\setminus L);R) @<<< H^i_c(U;R)\\
@. @A\psi^*_0AA @Vs^*VV\\
@. H^i_c(K_0\setminus L;R) @<=<< H^i_c(K_0\setminus L;R).\\
\endCD
$$
Therefore it suffices to prove Lemma for open sets $U$ of
the type $\psi_0^{-1}(K_0\setminus L)$ where $L$ is a subcomplex in $K_0$, and
for cohomology classes that came from
$H^i_c(K_0\setminus L;R)=H^i(K_0,L;R)$.

The sequence of covers $\{\sU_j\}$ defines an inverse sequence of simplicial
complexes
$$
K_0 @<\psi^1_0<< K_1 @<\psi^2_1<< K_2 @<\psi^3_2<<\dots
$$
with the limit space $X$ and with
$\psi_j=\psi^{j+1}_j\circ\psi_{j+1}$ for all $j$. We may assume
that every bonding map $\psi^{j+1}_j$ is simplicial with respect
to an iterated barycentric subdivision of the target $K_j$. If
$A=\lim_{\leftarrow}\{L_j,\psi^{j+1}_j|_{\dots}\}\subset_{Cl} X$,
then
$$
H^*_c(X\setminus A;R)=H^*(X,A;R)=\lim_{\rightarrow}H^*(K_j,L_j;R).
$$

By $C^i_j$  we denote the group of simplicial $R$-cochains on $K_j$
with compact supports. If $L\subset K_j$ is a subcomplex, by $C^i_j(K_j,L)$
we denote the subgroup of cochains which are zero on simplices from $L$.
For an $i$-simplex $\sigma\subset K_j$ we denote by $\phi_{\sigma}$ the cochain
that takes $\sigma$ to $1\in R$ and which is zero on all other simplices.
Then the cochains $\phi_{\sigma}$ form a basis of $C^i_j$.
Simplicial approximations $p^{j+1}_j:K_{j+1}\to K_j$ of the bonding maps
$\psi^{j+1}_j$ induce homomorphisms $(p^{j+1}_j)^*:C^*_j\to C^*_{j+1}$.
Let $A=\lim_{\leftarrow}\{L_j,\psi^{j+1}_j|_{\dots}\}$ be the limit
of inverse sequence of subcomplexes $L_j\subset K_j$. By the construction of
a simplicial approximation $p^{j+1}_j(L_{j+1})\subset L_j$.
The induced homomorphism on cohomology level
$$(p^{j+1}_j)^*:H^*(K_j,L_j;R)\to
H^*(K_{j+1},L_{j+1};R)$$ coincides with
$$(\psi^{j+1}_j)^*:H^*(K_j,L_j;R)\to
H^*(K_{j+1},L_{j+1};R)$$ because the maps
$$p^{j+1}_j,\psi^{j+1}_j:(K_{j+1},L_{j+1})\to (K_j,L_j)$$ are
homotopic.
Since the homology of a direct limit of chain complexes equals the direct
limit of homologies,
the limit group $H_c^i(X\setminus A;R)$ of the direct system
$$\{H^i(K_j,L_j;R), (p^{j+1}_j)^*\}$$ of the homology groups
of the cochain complexes $C^i_j(K_j,L_j)$ is the homology of the
cochain complex $\lim_{\rightarrow}C_k^i(K_j,L_j)$.
For very element $\phi\in\lim_{\rightarrow}C_j^i(K_j,L_j)$
we take a minimal $k=k_{\phi}$
such that $\phi$ comes from some $\phi_k\in C^i_k(K_k,L_k)$ and fix such
$\phi_k$. Note that if $p^{j+1}_j$ are surjections, then $\phi_k$ is unique.
We define
$supp_X\phi=\psi_k^{-1}(supp(\phi_k)).$

Let $\dim K_0=m$. Then the Lemma is trivial for $i>m$
(for our special setting).
So, if $m\le n$, Lemma holds true. Let $m>n$.
By Proposition 2 the group $H^i_c(X;R)$ is uniformly trivial for $n<i\le m$.
Let $\nu:\R_+\to\R_+$ be a corresponding function from for all $i$ between $n$
and $m$. Let $\rho(t)=\nu(2t)$. We use notation $\rho^k=\rho\circ\dots\circ\rho$
for the $k$th iteration. Let $p^*_j;C^*_j\to\lim_{\to}C^*_k$ denote the natural
homomorphism. 
Starting from $i=m$ by induction
for $i>n$ we construct a sequence of
homomorphisms
$$A_i:C^i_0\to\lim_{\rightarrow}C_k^{i-1}$$
such that
\roster
\item{} Being restricted to the subgroup of cocycles $Z^i_0\subset C^i_0$, the homomorphism
$A_i$ is a lift of
$p^*_0|_{Z^i_0}:Z^i_0\to\lim_{\to}C_k^i$ with respect to
$\delta:lim_{\rightarrow}C_k^{i-1}\to lim_{\rightarrow}C_k^{i}$;
\item{} $supp_X A_i(\phi_{\sigma})\subset \dot B_{\rho^{m-i+1}(d)}(x)$
where $\psi_0(x)\in\sigma$.
\endroster
Construction of $A_m$:
For every $m$-simplex $\sigma\subset K_0$ we consider the inverse sequence
$$
K_0\setminus\dot\sigma @<\psi^1_0<< (\psi^1_0)^{-1}(K_0\setminus\dot\sigma)
@<\psi^2_1<< (\psi^2_1)^{-1}(K_0\setminus\dot\sigma) @<<<\dots
$$
of simplicial complexes with the limit space $Y\subset X$.
Note that the set $X\setminus Y$
is contained
in an open $d$-ball $\dot B_{d}(x)$ with $\psi_0(x)\in\sigma$.
Since $H^m_c(X;R)$ is uniformly trivial with
the function
$\nu$, the inclusion $\dot B_{d}(x)\subset\dot B_{\nu(d)}(x)$
induces zero homomorphism. Thus, the cohomology class $[\phi_{\sigma}]\in
H^m(K_0,K_0\setminus\dot\sigma;R)$ defines the class
$(\psi_0^{-1})^*[\phi_{\sigma}]\in H^i_c(X\setminus Y;R)$ which
goes to zero under this inclusion $X\setminus Y\subset \dot B_{\nu(d)}(x)$.
Taking into account the direct limit description of the cohomologies 
we can say that
there is $j$ such that $(\psi_0^j)^*[\phi_{\sigma}]$ goes to zero
in $H^i(K_j,K_j\setminus st(\psi_j(B_{\nu(d)}(x)));R)$.
We can
rephrase this on the cochain level as follows:
There is a cochain $a\in C^{m-1}_j$
supported in  $\psi_j(B_{\nu(d)}(x))$
such that $\delta a=(p_0^j)^*(\phi_{\sigma})$.
We define $A_m(\phi_{\sigma})=p_j^*(a)$. Since $A_m$ is defined on
the basis and $C^m_0$ is a free $R$-module, it is defined everywhere.
Clearly, the conditions (1)-(2) hold.

Assume that $A_{i+1}$ is constructed.
Let $\sigma$ be an $i$-dimensional simplex.
We define $A_i$ on a basic cochain
$\phi_{\sigma}$ as follows. Note that
$b_{\sigma}=p_1^*(\phi_{\sigma})-A_{i+1}\delta\phi_{\sigma}$ is a cocycle in
$\lim_{\to}C_j^i$ and hence in defines a cocycle in $C_k^i$ for some $k$. 
For this cocycle we use the same notation $b_{\sigma}$.
By the condition (2),
$supp_X(b_{\sigma})\subset \dot B_{2\rho^{m-i}(d)}(x)$ where
$\psi_0(x)\in\sigma$. Hence
$b_{\sigma}$ defines a cohomology class in $H^i_c(\dot B_{2\rho^{m-i}(d)}(x);R)$. 
As above using uniform triviality
of the cohomology group $H^i_c(X;R)$ we obtain that for some $j>k$ we have
$\delta a=(p_k^j)^*(b_{\sigma})$ for some cochain $a\in C^{i-1}_j$ 
with the support in
$\psi_j(B_{\nu(2\rho^{m-i}(d))}(x))$. We define $A_i(\phi_{\sigma})=p_j^*(a)$.
It's easy to verify the conditions (1)-(2).

Let $r=2\rho^{m-n}(d)$.
If $L\subset K_0$ is a subcomplex and $\alpha\in C^i_0$, $i>n$, is
a cocycle that represents a class in $H^i(K_0,L;R)$, then $A_i(\alpha)$
comes from a cochain $\beta\in C^{i-1}_j$ for some $j$ that 
cobounds $(p_0^j)^*(\alpha)$.
By the condition (2) $\phi^{-1}_j(supp A_i(\alpha))\subset
N_r(\psi^{-1}_0(supp(\alpha)))$
\qed
\enddemo

\demo{Proof of Theorem 1}
The inequality $\dim_RZ\ge gcd_R(Z)$ is obvious.

We prove the reverse inequality. Let $\dim_RZ=n$
and assume that $gcd_R\Gamma < n$. Apply Lemma 1 and Lemma 2 to
obtain $z\in Z$, an open neighborhood
$V$ in $Z$ and an
open set $\bar U\subset\bar X$ with $\bar U\cap Z=V$ such that
the homomorphism $j:H_c^n(V;R)\to
H^{n+1}_c(U;R)$ is a nontrivial monomorphism.
We may assume  that $V$ is as in the definition of a dimensionally
essential point, i.e., for every smaller neighborhood $V'$,
$z\in V'\subset V$, the inclusion homomorphism
$H^n_c(V';R)\to H_c^n(V;R)$ is non zero.
Let $r$ be as in
Lemma 3. We take
a neighborhood $W$ of $z$ in $\bar X$ such that $Cl_{\bar X}W\subset \bar U$
and $N_{2r}(X\setminus U)\cap W=\emptyset$. Then by Axiom 4
the $r$-neighborhood $N_{r}(W\setminus Z)$ is contained in $U$.
We apply Lemma 2 to obtain an open  neighborhood of
$z$, $\bar U_0\subset W$,
such that $$j':H^n_c(V_0;R)\to H^{n+1}_c(U_0;R)$$ is a nonzero
monomorphism where
$V_0=\bar U_0\cap Z$ and $U_0=\bar U_0\setminus Z$.
There is an element $\alpha$ which goes to a non zero element under
the inclusion homomorphism $H^n_c(V_0;R)\to H^n_c(V;R)$.
Since $N_{r_0}(W\setminus Z)$ is contained in $U$, by Lemma 3 the
inclusion homomorphism
$H^{n+1}_c(U_0;R)\to H^{n+1}_c(U;R)$ is trivial. The commutative
diagram with injective vertical arrows
$$
\CD
H^n_c(V_0;R) @>>>  H^n_c(V;R)\\
@Vj'VV @VjVV\\
H^{n+1}_c(U_0;R) @>>>  H^{n+1}_c(U;R)\\
\endCD
$$
gives us a contradiction.
\qed
\enddemo

\proclaim{Corollary 2}
The cohomological dimension of the  boundary $\dim_R\partial\Gamma$ 
for groups with $\sZ$-structure is a quisi-isometry invariant.
\endproclaim
\demo{Proof}
Let $\Gamma_1$ and $\Gamma_2$ be two groups with $\sZ$-structures
$(\bar X_1, Z_1)$ and $(\bar X_2, Z_2)$ respectively.
By Theorem 1 and the cohomology exact sequence of pairs $(\bar X_i, Z_i)$
it suffices to show that $H^*_c(X_1;R)=H^*_c(X_2;R)$.
Since the spaces $X_i$, $i=1,2$, are uniformly contractible, by 
a theorem of Roe $H^*_c(X_i;R)=HX^*(X_i;R)$ [Roe], [HR] where $HX^*(X;R)$ denotes
the coarse cohomology groups of $X$. The chain of coarse equivalences
$$
X_1\sim\Gamma_1\sim\Gamma_2\sim X_2
$$
completes the proof.
\qed
\enddemo

\enddocument